*Review*

# A focus on the Riemann's hypothesis

### Jean-Max CORANSON-BEAUDU

Le Lamentin, Martinique, F.W.I, France.

**Riemann's hypothesis, formulated in 1859, concerns the location of the zeros of Riemann's Zeta function. The history of the Riemann hypothesis is well known. In 1859, the German mathematician B. Riemann presented a paper to the Berlin Academy of Mathematic. In that paper, he proposed that this function, called Riemann-zeta function takes values 0 on the complex plane when s=0.5+it. This hypothesis has great significance for the world of mathematics and physics. This solutions would lead to innumerable completions of theorems that rely upon its truth. Over a billion zeros of the function have been calculated by computers and shown that all are on this line s = 0.5+it. In this paper, we initially show that Riemann's ζ (Zêta) function and the analytical extension of this function called ℵ (Aleph)) are distinct. After extending this function in the complex plane except the point s=1, we will show the existence and then the uniqueness of real part zeros equal to 1/2.**

## INTRODUCTION

Riemann's hypothesis is expressed as following:

All non-trivial zeros of the function $\zeta(s)$ are located on the complex line $\Re(s) = \frac{1}{2}$

## INTRODUCTION - ON THE ANALYTICAL EXTENSION OF THE FUNCTION $\zeta$

The analytical extension of the function $\zeta(s)$ on $\mathbb{C}$ will be called $\aleph(s)$ in order to distinguish it from the function of Riemann. Riemann's Zeta function is written:

$$\zeta(s) = \sum_{n=1}^{\infty} \frac{1}{n^s}$$

For all complex numbers $\Re(s) > 1$

The function $\frac{1}{x^s}$ for $x \in \mathbb{R}$ and $s \in \mathbb{C}$ is differentiable p times. The p-th derivative of this function is written:

$$\left(\frac{1}{x^s}\right)^{'p} = (-1)^p s(s+1)(s+2)\ldots(s+p-1)\frac{1}{x^{s+p}}$$

Applying Euler Mac-Laurin's (Havil, 2003; Poels, 2011) formula to the function:

$$\sum_{n=1}^{N}\frac{1}{p^s} = \int_1^N \frac{dx}{x^s} + \frac{1}{2}\left(1 + \frac{1}{N^s}\right) - \sum_{j=1}^{M}\frac{b_{2j}}{(2j)!}s(s+1)\ldots(s+2j-2)\left(\frac{1}{N^{s+2j-1}} - 1\right) + R_M(s) \quad (1)$$

$$R_M(s) = -\frac{s(s+1)\ldots(s+2M-1)}{(2M+1)!}\int_1^N B^*_{2M+1}(x)x^{-s-2M-1}dx$$

*Corresponding author. E-mail: jmcoranson@me.com

Where $B^*_{2M+1}(x) = B_{2M+1}(x - E(x))$ is a 1-periodic function $B_n(x)$, called the p-th Bernoulli polynome and $b_n = B_n(0)$, called the p-th number of Bernoulli.

For $N \to +\infty$ the left member of the Equation 1 leans towards $\zeta(s)$ and the development of Euler MacLaurin, a right sided part of the equation is defined by:

$$\frac{1}{s-1} + \frac{1}{2} + \sum_{j=1}^{M} \frac{b_{2j}}{(2j)!} s(s+1)\ldots(s+2j-2) + \sigma_M(s)$$

With

$$\sigma_M(s) = -\frac{s(s+1)\ldots(s+2M)}{(2M+1)!} \int_1^{+\infty} B^*_{2M+1}(x) x^{-s-2M-1} dx$$

$\sigma_M$ being a convergent integral for $\Re(s) > 1 - 2M \; \forall M \in \mathbb{N}^*$, converges for all s of the complex plan except in $s = 1$.

The other members of MacLaurin's development being polynomes, the analytical extension of the Zeta function is defined by the entire complex plan except in 1. The analytical extension (Edwards, 1974; Lachaud, 2001) of Riemann's function is expressed by the following formula:

$$\aleph(s) = \begin{cases} \zeta(s) = \sum_{n=1}^{\infty} \frac{1}{n^s} & \text{for} \quad \Re(s) > 1 \\ \frac{1}{s-1} + \frac{1}{2} + \sum_{j=1}^{M} \frac{b_{2j}}{(2j)!} s(s+1)\ldots(s+2j-2) + \sigma_M(s) & \forall s \in \mathbb{C}/\{1\} \end{cases}$$

(2)

It is clear that calculating the value of $\zeta(s)$ for values such as 0 or -1 with the following formula,

$$\zeta(s) = \sum_{n=1}^{\infty} \frac{1}{n^s}$$

is impossible. So,

$$\zeta(-1) = \sum_{n=1}^{\infty} n = -\frac{1}{12}$$

is nonsense.

On the other $\aleph(-1)$ does exist through the converging integral $\sigma_M(-1)$.

The function $\zeta(s)$ does not admit zeros on its domain. $\Re(s) > 1$.

On the other hand $\aleph(s)$ being holomorphic on $\mathbb{C}/\{1\}$ there are zeroes for $\Re(s) \leq 1$.

**ON THE ZEROS OF THE FUNCTION $\aleph(s)$**

According to Fourier's (Andreas, 1987) analysis, the function $x \to e^{-\pi x^2}$ that belongs to Schwartz's (Schwartz, 1966) space of fast decay functions to infinity, coincides with his transformed Fourier, that is:

$$\int_{-\infty}^{+\infty} e^{-\pi x^2} e^{-2i\pi u x} \, dx = e^{-\pi u^2}$$

By making the variable change of $x \to \frac{x}{\sqrt{t}}$ in this integral, the Fourier transformation of the function

$$f(x) = e^{-\pi t x^2} \text{ is } \hat{f}(u) = \frac{1}{\sqrt{t}} e^{-\frac{\pi u^2}{t}}$$

And all functions of Schwartz's space we have the following relationship:

$$\forall n \in \mathbb{Z} \quad \sum_{n=-\infty}^{\infty} f(n) = \sum_{n=-\infty}^{\infty} \hat{f}(n)$$

which implies that

$$\forall t > 0 \; \mho(t) = \sum_{n=-\infty}^{\infty} e^{-\pi t n^2} = \frac{1}{\sqrt{t}} \sum_{n=-\infty}^{\infty} e^{-\frac{\pi n^2}{t}} \quad (3)$$

$\mho$ and $\psi$ functions meet the following functional equations

$$\forall u > 0 \; \mho(u) = \frac{1}{\sqrt{u}} \mho\left(\frac{1}{u}\right) \; et \; \psi(u) = \frac{\mho(u) - 1}{2} = \sum_{n=1}^{\infty} e^{-\pi u n^2}$$

$\psi$ checks:

$$\psi\left(\frac{1}{u}\right) = \frac{\mho\left(\frac{1}{u}\right) - 1}{2} = \frac{\mho(u)\sqrt{u} - 1}{2} = \frac{\sqrt{u}(2\psi(u) + 1) - 1}{2} = \sqrt{u}\psi(u) + \frac{\sqrt{u}}{2} - \frac{1}{2}$$

That is,

$$\forall u > 0 \; \psi(u) = \frac{1}{\sqrt{u}} \psi\left(\frac{1}{u}\right) + \frac{1}{2\sqrt{u}} - \frac{1}{2} \quad (4)$$

$\forall s \,/\, \Re(s) > 1, et \; n \neq 0$

To calculate the full one below by posing the variable change,

$$u = \frac{1}{\pi n^2} t$$

$$I_n = \int_0^\infty u^{\frac{s}{2}-1} e^{-\pi n^2 u} du = \int_0^\infty \frac{\pi n^2}{(\pi n^2)^{\frac{s}{2}+1}} t^{\frac{s}{2}-1} e^{-t} dt$$

$$= \frac{1}{\pi^{s/2}} \frac{1}{n^s} \int_0^\infty t^{\frac{s}{2}-1} e^{-t} dt$$

$$I_n = \frac{1}{\pi^{s/2}} \frac{1}{n^s} \Gamma\left(\frac{s}{2}\right)$$

By summing n, we obtain:

$$\sum_{n=1}^{n=\infty} I_n = \sum_{n=1}^\infty \int_0^\infty u^{\frac{s}{2}-1} e^{-\pi n^2 u} du = \sum_{n=1}^\infty \frac{1}{\pi^{s/2}} \frac{1}{n^s} \Gamma\left(\frac{s}{2}\right)$$

The inversion between infinite summation and integration is justified by the convergence properties of the function $e^{-\pi n^2 u}$. So we obtain:

$$\int_0^\infty u^{\frac{s}{2}-1} \sum_{n=1}^\infty e^{-\pi n^2 u} du = \pi^{-\frac{s}{2}} \Gamma\left(\frac{s}{2}\right) \sum_{n=1}^\infty \frac{1}{n^s} = \pi^{-\frac{s}{2}} \Gamma\left(\frac{s}{2}\right) \zeta(s)$$

That is,

$$\int_0^\infty u^{\frac{s}{2}-1} \psi(u) du = \pi^{-\frac{s}{2}} \Gamma\left(\frac{s}{2}\right) \zeta(s) \quad (5)$$

With $\zeta(s)$ is the function of Riemann for $\Re(s) > 1$

The integral 5 is developed on the intervals, $[0;1] \cup [1;+\infty]$. We have:

$$\pi^{-\frac{s}{2}} \Gamma\left(\frac{s}{2}\right) \zeta(s) = \int_0^\infty u^{\frac{s}{2}-1} \psi(u) du = \int_0^1 u^{\frac{s}{2}-1} \psi(u) du + \int_1^\infty u^{\frac{s}{2}-1} \psi(u) du$$

as

$$\psi(u) = \frac{1}{\sqrt{u}} \psi\left(\frac{1}{u}\right) + \frac{1}{2\sqrt{u}} - \frac{1}{2}$$

So on the interval $[0;1]$ we can write:

$$\int_0^1 u^{\frac{s}{2}-1} \psi(u) du = \int_0^1 u^{\frac{s}{2}-1} \left(\frac{1}{\sqrt{u}} \psi\left(\frac{1}{u}\right) + \frac{1}{2\sqrt{u}} - \frac{1}{2}\right) du$$

By placing $u = \frac{1}{v}$ in the first part of the integral we have:

$$\int_0^1 u^{\frac{s}{2}-1} \psi(u) du = -\int_\infty^1 v^{\frac{-s}{2}+1} \left(v^{\frac{1}{2}} \psi(v)\right) \frac{1}{v^2} dv + \int_0^1 u^{\frac{s}{2}-1} \left(-\frac{1}{2} + \frac{1}{2\sqrt{u}}\right) du$$

$$\int_0^1 u^{\frac{s}{2}-1} \psi(u) du = \int_1^\infty \psi(v) v^{\frac{-s}{2}+1-2+\frac{1}{2}} dv - \frac{u^{\frac{s}{2}}}{2\left(\frac{s}{2}\right)}\Bigg|_0^1 + \frac{u^{\frac{s}{2}-\frac{1}{2}}}{2\left(\frac{s-1}{2}\right)}\Bigg|_0^1$$

$$\int_0^1 u^{\frac{s}{2}-1} \psi(u) du = \int_1^\infty \psi(u) u^{\frac{-s}{2}-\frac{1}{2}} du - \frac{1}{s} + \frac{1}{s-1}$$

Therefore

$$\pi^{-\frac{s}{2}} \Gamma\left(\frac{s}{2}\right) \zeta(s) = \int_1^\infty \left(u^{\frac{s}{2}-1} + u^{\frac{-s-1}{2}}\right) \psi(u) du - \frac{1}{s} - \frac{1}{1-s} \quad (6)$$

This integral is converging for any complex except 0 and 1. $\zeta(s)$ function is defined by continuity on $\mathbb{C}/\{0;1\}$ as

$$\frac{1}{s} + \frac{1}{1-s} = \frac{1}{s(1-s)}$$

by multiplying Equation 6 by $s(s-1)$ we have:

$$\pi^{-\frac{s}{2}} \Gamma\left(\frac{s}{2}\right) \zeta(s) s(s-1) = s(s-1) \int_1^\infty \left(u^{\frac{s}{2}-1} + u^{\frac{-s-1}{2}}\right) \psi(u) du + 1$$

therefore we use the term $\aleph(s)$ instead of $\zeta(s)$ and define:

$$\beth(s) = \pi^{-\frac{s}{2}} \Gamma\left(\frac{s}{2}\right) \aleph(s) s(s-1) = s(s-1) \int_1^\infty \left(u^{\frac{s}{2}-1} + u^{\frac{-s-1}{2}}\right) \psi(u) du + 1$$

as

$$\frac{s}{2} \Gamma\left(\frac{s}{2}\right) = \Gamma\left(\frac{s}{2}+1\right)$$

Then,

$$\beth(s) = 2\pi^{-\frac{s}{2}} \Gamma\left(\frac{s}{2}+1\right) \aleph(s)(s-1) = s(s-1) \int_1^\infty \left(u^{\frac{s}{2}-1} + u^{\frac{-s-1}{2}}\right) \psi(u) du + 1 \quad (7)$$

This integral is defined $\forall s \in \mathbb{C}$ thanks to the rapid decay property of the $\psi$ function to infinity. It can be said that $\beth(s)$ is holomorphic in $\mathbb{C}$.

So, $\beth(s) = \Phi(s) \aleph(s)$, with meromorphic $\Phi(s) = 2\pi^{-\frac{s}{2}} \Gamma\left(\frac{s}{2}+1\right)(s-1)$ and $\beth(s)$ is holomorphic, then $\aleph(s)$ is meromorphic.

On the other hand, $\beth(s) = \beth(1-s)$ is a functional relationship between $\aleph(s)$ and $\aleph(1-s)$:

$$\pi^{-\frac{s}{2}}\Gamma\left(\frac{s}{2}\right)\aleph(s)s(s-1) = \pi^{-\frac{1-s}{2}}\Gamma\left(\frac{1-s}{2}\right)\aleph(1-s)s(s-1)$$

That is,

$$\Gamma\left(\frac{S}{2}\right)\aleph(s) = \pi^{s-\frac{1}{2}}\Gamma\left(\frac{1-s}{2}\right)\aleph(1-s) \quad (8)$$

The function $\beth$ is written on $\mathbb{C}$

$$\beth(s) = s(s-1)\int_1^\infty \left(u^{\frac{s-\frac{1}{2}}{2}} + u^{-\frac{s-\frac{1}{2}}{2}}\right)u^{-\frac{3}{4}}\psi(u)du + 1$$

That is,

$$\beth(s) = 2\int_1^\infty s(s-1)u^{-\frac{3}{4}}\psi(u)\cosh\left[\left(s-\frac{1}{2}\right)\frac{\ln(u)}{2}\right]du + 1 \quad (9)$$

verify that,

$$\beth(s) = \beth(1-s) \text{ and } \beth(0) = \beth(1-0) = 1$$

**The trivial zeroes**

$$\beth(s) = 2\pi^{-\frac{s}{2}}\Gamma\left(\frac{S}{2}+1\right)\aleph(s)(s-1) \Rightarrow \aleph(s) = \frac{1}{\Gamma\left(\frac{S}{2}+1\right)}\frac{\pi^{\frac{s}{2}}}{2(s-1)}\beth(s)$$

The function $\frac{1}{\Gamma\left(\frac{s}{2}+1\right)}$ included as zeroes $\frac{s}{2}+1 = -k$, $k \in \mathbb{N}$; $s = -2(k+1)$. On $\mathbb{C}/\{1\}$

$\frac{\pi^{\frac{s}{2}}}{2(s-1)}\beth(s)$ function is holomorphic.

Therefore, the function $\aleph(s)$ included the same trivial zeroes as the zeroes of function $\frac{1}{\Gamma\left(\frac{s}{2}+1\right)}$

$s = -2(1+k), k \in \mathbb{N}$ which are whole negative pairs.

**Non-trivial zeroes**

If there are non-trivial zeroes in the complex plan for this function $\beth$,

We expressed them as $z_k = a_k + ib_k$ $k \in \mathbb{N}$ and these are the same zeroes as the function $\aleph(s)$. Note $\Re(.)$ the real part and $\Im(.)$ the imaginary part.

These zeroes check the next relationship for the $\beth$ function.

$$\beth(z_k) = 0 \Leftrightarrow \begin{cases} \Re(\beth(z_k)) = 0 \\ \Im(\beth(z_k)) = 0 \end{cases} \forall k \in \mathbb{N} \quad (10)$$

By writing the real and imaginary part of the integral 9

for $s = z_k$ we obtain:

$$\begin{cases} 2\int_1^\infty u^{-\frac{3}{4}}\psi(u)\Re\left[z_k(z_k-1)\cosh\left[\left(z_k-\frac{1}{2}\right)\frac{\ln(u)}{2}\right]\right]du + 1 = 0 \\ 2\int_1^\infty u^{-\frac{3}{4}}\psi(u)\Im\left[z_k(z_k-1)\cosh\left[\left(z_k-\frac{1}{2}\right)\frac{\ln(u)}{2}\right]\right]du = 0 \end{cases} \quad (11)$$

We seek to identify complex $z_k$ values that verify the Equation 11 as,

$$z_k(z_k - 1) = (a_k + ib_k)(a_k - 1 + ib_k) = a_k(a_k - 1) - b_k^2 + ib_k(2a_k - 1)$$

and,

$$\cosh\left[\left(z_k - \frac{1}{2}\right)\frac{\ln(u)}{2}\right] = \cosh\left[\left(a_k - \frac{1}{2} + ib_k\right)\frac{\ln(u)}{2}\right]$$
$$= \cosh\left[\left(a_k - \frac{1}{2}\right)\frac{\ln(u)}{2}\right]\cos\left(b_k\frac{\ln(u)}{2}\right) + i\sinh\left[\left(a_k - \frac{1}{2}\right)\frac{\ln(u)}{2}\right]\sin\left(b_k\frac{\ln(u)}{2}\right)$$

We note $R(a_k, b_k, u)$ the real part of the product $z_k(z_k - 1)\cosh\left[\left(z_k - \frac{1}{2}\right)\frac{\ln(u)}{2}\right]$

And $I(a_k, b_k, u)$ the imaginary part of the product $z_k(z_k - 1)\cosh\left[\left(z_k - \frac{1}{2}\right)\frac{\ln(u)}{2}\right]$

We have got:

$$R(a_k, b_k, u) = (a_k(a_k - 1) - b_k^2)\cosh\left[\left(a_k - \frac{1}{2}\right)\frac{\ln(u)}{2}\right]\cos\left(b_k\frac{\ln(u)}{2}\right)$$
$$- b_k(2a_k - 1)\sinh\left[\left(a_k - \frac{1}{2}\right)\frac{\ln(u)}{2}\right]\sin\left(b_k\frac{\ln(u)}{2}\right)$$

From this expression, we obtain, because of property of $A\cos\left[b_k\frac{\ln(u)}{2}\right] + B\sinh\left[b_k\frac{\ln(u)}{2}\right]$,

$$R(a_k, b_k, u) = \sqrt{A^2 + B^2}\sin\left[b_k\frac{\ln(u)}{2} + \frac{\pi}{2}Sign(A) - \arctan\left(\frac{B}{A}\right)\right] \quad (12)$$

With
$A = (a_k(a_k - 1) - b_k^2)\cosh\left[\left(a_k - \frac{1}{2}\right)\frac{\ln(u)}{2}\right]$ and $B = -b_k(2a_k - 1)\sinh\left[\left(a_k - \frac{1}{2}\right)\frac{\ln(u)}{2}\right]$

We also have:

$$I(a_k, b_k, u) = (a_k(a_k - 1) - b_k^2)\sinh\left[\left(a_k - \frac{1}{2}\right)\frac{\ln(u)}{2}\right]\sin\left(b_k\frac{\ln(u)}{2}\right)$$
$$+ b_k(2a_k - 1)\cosh\left[\left(a_k - \frac{1}{2}\right)\frac{\ln(u)}{2}\right]\cos\left(b_k\frac{\ln(u)}{2}\right)$$

$$I(a_k, b_k, u) = \sqrt{U^2 + V^2}\sin\left[b_k\frac{\ln(u)}{2} + \frac{\pi}{2}Sign(U) - \arctan\left(\frac{V}{U}\right)\right] \quad (13)$$

With $U = b_k(2a_k - 1)\cosh\left[\left(a_k - \frac{1}{2}\right)\frac{\ln(u)}{2}\right]$ and $V = (a_k(a_k - 1) - b_k^2)\sinh\left[\left(a_k - \frac{1}{2}\right)\frac{\ln(u)}{2}\right]$

For the imaginary part of the integral of the equation

system 11 we have:

$$2\int_1^\infty u^{-\frac{3}{4}}\psi(u)\Im\left[z_k(z_k-1)\cosh\left[\left(z_k-\frac{1}{2}\right)\frac{\ln(u)}{2}\right]\right]du = 0$$

That is,

$$2\int_1^\infty u^{-\frac{3}{4}}\psi(u)I(a_k,b_k,u)\,du = 0 \qquad (14)$$

Are there couples $(a_k, b_k)$ such as Equation 14 equals zero?

Because of the convergence characteristics of the integral, we have the following property:

$$-\int_1^\infty u^{-\frac{3}{4}}\psi(u)|I(a_k,b_k,u)|\,du \le \int_1^\infty u^{-\frac{3}{4}}\psi(u)I(a_k,b_k,u)\,du \le \int_1^\infty u^{-\frac{3}{4}}\psi(u)|I(a_k,b_k,u)|\,du$$

That is,

$$-\int_1^\infty u^{-\frac{3}{4}}\psi(u)\sqrt{U^2+V^2}\,du \le \int_1^\infty u^{-\frac{3}{4}}\psi(u)I(a_k,b_k,u)\,du \le \int_1^\infty u^{-\frac{3}{4}}\psi(u)\sqrt{U^2+V^2}\,du$$

Applying the properties of the full continuous and positive function, and the squeeze theorem,

we have:

$$\int_1^\infty u^{-\frac{3}{4}}\psi(u)\sqrt{U^2+V^2}\,du = 0 \Rightarrow \sqrt{U^2+V^2} = 0 \Leftrightarrow \begin{cases} U = 0 \\ V = 0 \end{cases} \forall\, u \ge 1$$

The existence of the couples $(a_k, b_k)$ such as:

$$\begin{cases} U = 0 \Leftrightarrow a_k = \frac{1}{2}\ ou\ b_k = 0 \\ V = 0 \Leftrightarrow a_k = \frac{1}{2}\ ou\ a_k(a_k-1) - b_k^2 = 0 \end{cases}$$

The system is reduced to three pairs of solutions:

$$\begin{cases} a_k = \frac{1}{2} \\ b_k \in \mathbb{R} \end{cases} \text{or} \begin{cases} a_k = 1 \\ b_k = 0 \end{cases} \text{or} \begin{cases} a_k = 0 \\ b_k = 0 \end{cases}$$

We're checking that:

$z_0 = (0,0)\ ou\ z_1 = (1,0)$ are des trivial solutions $\Im(\beth(0)) = \Im(\beth(1)) = 0$ because $\beth(0) = \beth(1) = 1$

And $\begin{cases} a_k = \frac{1}{2} \\ b_k \in \mathbb{R} \end{cases}$ are non-trivial zeroes. As a result, we have shown that there are non-trivial zeroes on the critical axis $\Re(s) = \frac{1}{2}$.

The imaginary part $b_k$ of these zeroes is identified using the first integral of the equation system 11 expressed:

$$2\int_1^\infty u^{-\frac{3}{4}}\psi(u)R(a_k,b_k,u)\,du + 1 = 0 \qquad (15)$$

That is,

$$2\int_1^\infty u^{-\frac{3}{4}}\psi(u)\sqrt{A^2+B^2}\sin\left[b_k\frac{\ln(u)}{2} + \frac{\pi}{2}Sign(A) - \arctan\left(\frac{B}{A}\right)\right]du + 1 = 0$$

Taking into consideration the result found for the imaginary part of the integral, the following couples:

$$\begin{cases} a_k = \frac{1}{2} \\ b_k \in \mathbb{R} \end{cases}$$

We've got:

$$\begin{cases} A = (a_k(a_k-1) - b_k^2)\cosh\left[\left(a_k - \frac{1}{2}\right)\frac{\ln(u)}{2}\right] = -\frac{1}{4} - b_k^2 \\ B = b_k(2a_k - 1)\sinh\left[\left(a_k - \frac{1}{2}\right)\frac{\ln(u)}{2}\right] = 0 \end{cases}$$

Therefore the Equation 15 is written:

$$2\int_1^\infty u^{-\frac{3}{4}}\psi(u)\left(\frac{1}{4} + b_k^2\right)\cos\left(b_k\frac{\ln(u)}{2}\right)du = 1 \qquad (16)$$

and $b_k$ is a solution of the Equation 16. So there is an infinity of zeroes on the critical axis of the $\Re(s) = 1/2$ for $z_k = \frac{1}{2} + ib_k$

We show that these zeroes are all on the critical axis

Assumptions: *Suppose there are zeroes outside the critical axis and in the critical band.*

These zeros are written from the existing zeros on the critical line: with $y_k = z_k + \varepsilon e^{i\delta}$

$0 < \varepsilon < \frac{1}{2}$

We know that all $\overline{z_k} = 1 - z_k$ because $\Re(z_k) = \frac{1}{2}$   $k \in \mathbb{N}$

$\beth(z)$ is holomorph in $\mathbb{C}$, thus being an entire function. Weierstrass's factorization theorem (Patterson, 1995; Vento, 2003) states that any entire function can be represented by an infinite polynomial product with its zeroes. There is $g$ holomorph in $\mathbb{C}$ that does not cancel in $z_k$ and $\overline{z_k}$ such as:

$$\beth(z) = \prod_{k=1}^\infty \left(1 - \frac{z}{z_k}\right)\left(1 - \frac{z}{\overline{z_k}}\right)g(z(1-z)) = \prod_{k=1}^\infty \left(1 - \frac{z}{|z_k|^2}(1-z)\right)g(z(1-z)) \qquad (17)$$

We check that the $z_k$ and $\overline{z_k}$ are zeros of $\beth(z)$, the function verifies

$$\beth(1-z) = \beth(z)$$

$$\overline{\beth(z)} = \beth(\overline{z})$$

Suppose that $y_k = z_k + \varepsilon_k e^{i\delta_k}$ and $\overline{y_k} = \overline{z_k} + \varepsilon_k e^{-i\delta_k}$ are also zeros of $\beth(z)$, so we have:

$$\beth(z) = \beth_{z_k}(z) \prod_{k=1}^{\infty}\left(1 - \frac{z}{|y_k|^2}(y_k + \overline{y_k} - z)\right) g^*(z(1-z)) \quad (18)$$

with

$$\beth_{z_k}(z) = \prod_{k=1}^{\infty}\left(1 - \frac{z}{|z_k|^2}(1-z)\right)$$

And $g^*$ is holomorphic and does not cancel out for $y_k$, $z_k$, and their conjugates.

$$\beth(1-z) = \beth_{z_k}(1-z) \prod_{k=1}^{\infty}\left(1 - \frac{(1-z)}{|y_k|^2}(y_k + \overline{y_k} - (1-z))\right) g^*(z(1-z))$$

As $\beth_{z_k}(1-z) = \beth_{z_k}(z)$ then:

$$\beth(1-z) = \beth_{z_k}(z) \prod_{k=1}^{\infty}\left(1 - \frac{y_k + \overline{y_k} - 1}{|y_k|^2} - \frac{z}{|y_k|^2}(1 - y_k - \overline{y_k} + 1 - z)\right) g^*(z(1-z))$$

And $\beth(1-z) = \beth(z) \Leftrightarrow y_k + \overline{y_k} - 1 = 0$
So $z_k + \varepsilon_k e^{i\delta_k} + \overline{z_k} + \varepsilon_k e^{-i\delta_k} - 1 = 0$ that is,
$\varepsilon_k e^{i\delta_k} + \varepsilon_k e^{-i\delta_k} = 0$

Which is impossible since $\varepsilon_k \neq 0$

Therefore the hypothesis of zeros outside the critical axis leads to a contradiction in relation to the symmetries of function $\beth(X)$ in the critical band.

There are no zeroes outside the axis $\Re(z_k) = 1/2$.

## Conclusion

We have demonstrated:

(i) that the holomorphic function $\beth(s)$ had the same zeros as the function $\aleph(s)$ which is an analytical extension of Riemann's $\zeta(s)$ function because $\aleph(s) = \frac{1}{\Gamma(\frac{s}{2}+1)} \frac{\pi^{\frac{s}{2}}}{2(s-1)} \beth(s)$.

This result well known by the mathematical world, served us to find a holomorphic function simpler to exploit at the roots.
(ii) using the squeeze theorem on the integral form of the Riemann function, we show that there are a pairs $(a_k, b_k)$ that are zeros of the Riemann function and these zeros are on the line $s = \frac{1}{2} + it$

(iii) as Hadamard (1896) Charles-Jean (1916) have each proved that no zero of the analytical extension of the Zeta function could be found on the line Re(s)= 1, and therefore that all non-trivial zeroes must be in the interior of the critical band.
(iv) we have been hypothesis that if there were zeros, $y_k = z_k + \varepsilon_k e^{i\delta_k}$, in the critical band, with $0 < \varepsilon_k < \frac{1}{2}$, then this hypothesis leads to a contradiction. We used the Weierstrass's factorization theorem of holomorphic functions for $\beth(s)$, and applying functional relationship of symmetry, $\beth(1-z) = \beth(z)$, to demonstrate contradiction. Therefore, all non-trivial zeroes of $\beth$ are non-trivial zeroes of the analytical extension of the function $\zeta$ and have a real part $\frac{1}{2}$. These zeroes, noted $z_k = a_k + ib_k$ check the equation systems below:

$$\forall k \in \mathbb{Z} \begin{cases} a_k = \frac{1}{2} \\ b_k \in \mathbb{R} \,/\, 2\left(\frac{1}{4} + b_k^2\right) \int_1^\infty u^{-\frac{3}{4}} \psi(u) \cos\left[2b_k \frac{\ln(u)}{2}\right] du = 1 \end{cases}$$
(19)

## A simple digital example

A numerical integration by Rombert's method with order precision 5 and 20 iterations, we find the results of the complete system 16 with an error of $10^{-6}$.

$b_1 = 14.13472$
$b_2 = 21.02203$
$b_3 = 25.01085$
$b_4 = 30.42487$
$b_5 = 32.93506$

## CONFLICT OF INTERESTS

The author has not declared any conflict of interests.